\documentclass[a4paper]{article}
\usepackage{latexsym}

\title{Two distance-regular graphs}
\author{Andries E. Brouwer \& Dmitrii V. Pasechnik}

\newtheorem{Theorem} {Theorem} [section]
\newtheorem{Proposition} [Theorem] {Proposition}
\newcommand{\Proof}{ \noindent{\bf Proof:}\quad }
\newcommand{\qed}{\hfill$\Box$\medskip}

\newcommand{\Ff}{{\bf F}}
\newcommand{\F}{{\cal F}}

\newcommand{\adj}{\sim}
\newcommand{\mi}{{\hbox{$-$}}}
\newcommand{\pl}{{\hbox{$+$}}}
\date{2011/06/11}
\begin{document}
\maketitle
\begin{abstract}
We construct two families of distance-regular graphs,
namely the subgraph of the dual polar graph of type $B_3(q)$
induced on the vertices far from a fixed point,
and the subgraph of the dual polar graph of type $D_4(q)$
induced on the vertices far from a fixed edge.
The latter is the extended bipartite double of the former.
\end{abstract}

\section{The extended bipartite double}
We shall use $\adj$ to indicate adjacency in a graph.
For notation and definitions of concepts related to
distance-regular graphs, see \cite{BCN}.
We repeat the definition of extended bipartite double.

\medskip
The {\em bipartite double} of a graph $\Gamma$ with vertex set $X$
is the graph with vertex set $\{ x^+,x^- \mid x \in X \}$
and adjacencies $x^\delta \adj y^\epsilon$
iff $\delta\epsilon = -1$ and $x \adj y$. The bipartite double of a graph
$\Gamma$ is bipartite, and it is connected iff $\Gamma$ is connected and
not bipartite. If $\Gamma$ has spectrum $\Phi$, then its bipartite double
has spectrum $(-\Phi) \cup \Phi$. See also \cite{BCN}, Theorem 1.11.1.

\medskip
The {\em extended bipartite double} of a graph $\Gamma$ with vertex set $X$
is the graph with vertex set $\{ x^+,x^- \mid x \in X \}$, and the same
adjacencies as the bipartite double, except
that also $x^- \adj x^+$ for all $x \in X$.
The extended bipartite double of a graph $\Gamma$ is bipartite,
and it is connected iff $\Gamma$ is connected.
If $\Gamma$ has spectrum $\Phi$, then its extended bipartite double
has spectrum $(-\Phi-1) \cup (\Phi+1)$. See also \cite{BCN}, Theorem 1.11.2.

\section{Far from an edge in the dual polar graph of type $D_4(q)$}\label{bip}
Let $V$ be a vector space of dimension 8 over a field $F$,
provided with a nondegenerate quadratic form of maximal Witt index.
The maximal totally isotropic subspaces of $V$ (of dimension 4)
fall into two families $\F_1$ and $\F_2$, where the dimension of
the intersection of two elements of the same family is even (4 or 2 or 0)
and the dimension of the intersection of two elements of different families
is odd (3 or 1).

The geometry of the totally isotropic subspaces of $V$, where
$A \in \F_1$ and $B \in \F_2$ are incident when $\dim A \cap B = 3$
and otherwise incidence is symmetrized inclusion, is known as
the geometry $D_4(F)$.
The bipartite incidence graph on the maximal totally isotropic subspaces
is known as the dual polar graph of type $D_4(F)$.

Below we take $F = \Ff_q$, the finite field with $q$ elements, so that
graph and geometry are finite. We shall use projective terminology,
so that 1-spaces, 2-spaces and 3-spaces are called points, lines and planes.
Two subspaces are called disjoint when they have no point in common, i.e.,
when the intersection has dimension 0.

\begin{Proposition}\label{D4prop}
Let $\Gamma$ be the dual polar graph of type $D_4(\Ff_q)$.
Fix elements $A_0 \in \F_1$ and $B_0 \in \F_2$ with $A_0 \adj B_0$.
Let $\Delta$ be the subgraph of $\Gamma$ induced on the set of
vertices disjoint from $A_0$ or $B_0$. Then $\Delta$ is
distance-regular with intersection array
$\{q^3, q^3-1, q^3-q, q^3-q^2+1;~ 1,q,q^2-1,q^3\}$.
\end{Proposition}

The distance distribution diagram is

{\footnotesize

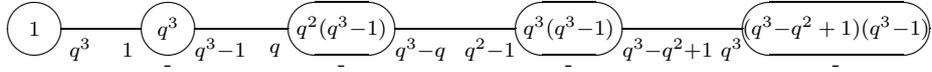
$$
\begin{picture}(345,20)(0,-16)
\put(10,0){\circle{20}}
\put(10,0){\makebox(0,0){1}}
\put(20,0){\line(1,0){30}}
\put(23,-8){\makebox(0,0)[l]{$q^3$}}
\put(47,-8){\makebox(0,0)[r]{1}}
\put(60,0){\circle{20}}
\put(60,0){\makebox(0,0){$q^3$}}
\put(60,-14){\makebox(0,0){-}}
\put(70,0){\line(1,0){35}}
\put(70,-8){\makebox(0,0)[l]{$q^3\mi1$}}
\put(102,-8){\makebox(0,0)[r]{$q$}}
\put(125,0){\oval(40,20)}
\put(125,0){\makebox(0,0){$q^2(q^3\mi1)$}}
\put(125,-14){\makebox(0,0){-}}
\put(145,0){\line(1,0){45}}
\put(145,-8){\makebox(0,0)[l]{$q^3\mi q$}}
\put(190,-8){\makebox(0,0)[r]{$q^2\mi1$}}
\put(210,0){\oval(40,20)}
\put(210,0){\makebox(0,0){$q^3(q^3\mi1)$}}
\put(210,-14){\makebox(0,0){-}}
\put(230,0){\line(1,0){45}}
\put(230,-8){\makebox(0,0)[l]{$q^3\mi q^2\pl1$}}
\put(275,-8){\makebox(0,0)[r]{$q^3$}}
\put(310,0){\oval(70,20)}
\put(310,0){\makebox(0,0){$(q^3\mi q^2+1)(q^3\mi1)$}}
\put(310,-14){\makebox(0,0){-}}
\end{picture}
$$
\vskip -1mm}
\Proof
There are $q^6$ elements $A \in  \F_1$ disjoint from $A_0$
and the same number of $B \in \F_2$ disjoint from $B_0$,
so that $\Delta$ has $2q^6$ vertices.

Given $A \in \F_1$, there are $q^3+q^2+q+1$ elements $B \in \F_2$
incident to it.
Of these, $q^2+q+1$ contain the point $A \cap B_0$ and hence are not
vertices of $\Delta$. So, $\Delta$ has valency $q^3$.

Two vertices $A,A' \in \F_1$ have distance 2 in $\Delta$ if and only if
they meet in a line, and the line $L = A \cap A'$ is disjoint
from $B_0$. If this is the case, then $L$ is in $q+1$ elements $B \in \F_2$,
one of which meets $B_0$, so that $A$ and $A'$ have $c_2=q$ common
neighbours in $\Delta$.

Given vertices $A \in \F_1$ and $B \in \F_2$ that are nonadjacent, i.e.,
that meet in a single point $P$, the neighbours $A'$ of $B$ at distance 2
to $A$ in $\Delta$ correspond to the lines $L$ on $P$ in $A$ disjoint
from $B_0$ and nonorthogonal to the point $A_0 \cap B$.
There are $q^2+q+1$ lines $L$ on $P$ in $A$, $q+1$ of which are
orthogonal to the point $A_0 \cap B$, and one further of which meets $B_0$.
(Note that the points $A_0 \cap B$ and $A \cap B_0$ are nonorthogonal
since neither point is in the plane $A_0 \cap B_0$ and $V$ does not
contain totally isotropic 5-spaces.)
It follows that $c_3 = q^2-1$, and also that $\Delta$ has diameter 4,
and is distance-regular.
\qed
%

\medskip
The geometry induced by the incidence relation of $D_4(F)$ on the vertices
of $\Delta$, together with the points and lines contained in the planes
disjoint from $A_0 \cup B_0$,
has Buekenhout-Tits diagram (cf.~\cite{P})

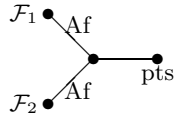
\begin{center}
\begin{picture}(60,30)(0,-15)
\put(0,0){\circle*{3.5}}
\put(-24,0){\circle*{3.5}}
\put(-41,17){\circle*{3.5}}
\put(-41,-17){\circle*{3.5}}
\put(0,0){\line(-1,0){24}}
\put(-24,0){\line(-1,1){17}}
\put(-24,0){\line(-1,-1){17}}
\put(-50,17){\makebox(0,0)[c]{\small $\F_1$}}
\put(-50,-17){\makebox(0,0)[c]{\small $\F_2$}}
\put(0,-6){\makebox(0,0)[c]{\small pts}}
\put(-30,13){\makebox(0,0)[c]{\small Af}}
\put(-30,-12){\makebox(0,0)[c]{\small Af}}
\end{picture}
\end{center}

\noindent
that is, the residue of an object $A \in \F_1$ is an affine 3-space,
where the objects incident to $A$ in $\F_2$ play the r\^{o}le of points.
Similar things hold more generally for $D_n(F)$ with arbitrary $n$,
and even more generally for all diagrams of spherical type.
See also \cite{BB}, Theorem 6.1.

\medskip
Let $P$ be a nonsingular point, and let $\phi$ be the reflection in the
hyperplane $H = P^\perp$. Then $\phi$ is an element of order two
of the orthogonal group that fixes $H$ pointwise, and consequently
interchanges $\F_1$ and $\F_2$. For each $A \in \F_1$ we have
$\phi(A) \adj A$. The quotient $\Gamma/\phi$ is the dual polar graph
of type $B_3(q)$, and we see that more generally the dual polar graph
of type $D_{m+1}(q)$ is the extended bipartite double of the dual polar
graph of type $B_m(q)$. 
The quotient $\Delta/\phi$ is a new distance-regular graph discussed
in the next section. It is the subgraph consisting of the vertices
at maximal distance from a given point in the dual polar graph
of type $B_3(q)$.
For even $q$ we have $B_3(q) = C_3(q)$, and it follows that the
symmetric bilinear forms graph on $\Ff_q^3$ is distance-regular,
see \cite{BCN} Proposition 9.5.10 and the diagram there on p. 286.

\section{Far from a point in the dual polar graph of type $B_3(q)$}\label{prim}

First a very explicit version of the graph of this section.

\begin{Proposition}\label{B3expl}
(i)
Let $W$ be a vector space of dimension $3$ over the field $\Ff_q$,
provided with an outer product $\times$.
Let $Z$ be the graph with vertex set $W \times W$ where
$(u,u') \adj (v,v')$ if and only if
$(u,u') \ne (v,v')$ and $u \times v + u' - v' = 0$.
Then $Z$ is distance-regular of diameter $3$ on $q^6$ vertices.
It has intersection array
$\{ q^3 - 1, q^3 - q, q^3 - q^2 + 1 ; 1, q, q^2 - 1 \}$
and eigenvalues
$q^3-1$, $q^2-1$, $-1$, $-q^2-1$
with multiplicities $1$, $\frac{1}{2}q(q+1)(q^3-1)$,
$(q^3-q^2+1)(q^3-1)$, $\frac{1}{2}q(q-1)(q^3-1)$, respectively.

(ii)
The extended bipartite double $\hat{Z}$ of $Z$
is distance-regular
with intersection array
$\{q^3, q^3-1, q^3-q, q^3-q^2+1;~ 1,q,q^2-1,q^3\}$
and eigenvalues $\pm q^3$, $\pm q^2$, $0$
with multiplicities $1$, $q^2(q^3-1)$, $2(q^3-q^2+1)(q^3-1)$,
respectively.

(iii)
The distance-1-or-2 graph $Z_1 \cup Z_2$ of $Z$,
which is the halved graph of $\hat{Z}$, is strongly regular
with parameters
$(v,k,\lambda,\mu) = (q^6,q^2(q^3-1),q^2(q^2+q-3),q^2(q^2-1))$.

\end{Proposition}

The distance distribution diagram of $Z$ is

{\footnotesize

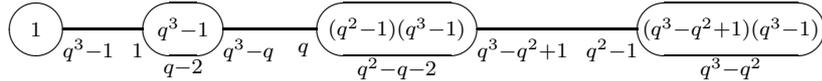
$$
\begin{picture}(345,20)(-10,-16)
\put(10,0){\circle{20}}
\put(10,0){\makebox(0,0){1}}
\put(20,0){\line(1,0){30}}
\put(20,-8){\makebox(0,0)[l]{$q^3\mi1$}}
\put(50,-8){\makebox(0,0)[r]{1}}
\put(65,0){\oval(30,20)}
\put(65,0){\makebox(0,0){$q^3\mi1$}}
\put(65,-15){\makebox(0,0){$q\mi2$}}
\put(80,0){\line(1,0){35}}
\put(80,-8){\makebox(0,0)[l]{$q^3\mi q$}}
\put(112,-8){\makebox(0,0)[r]{$q$}}
\put(145,0){\oval(60,20)}
\put(145,0){\makebox(0,0){$(q^2\mi1)(q^3\mi1)$}}
\put(145,-15){\makebox(0,0){$q^2\mi q\mi2$}}
\put(175,0){\line(1,0){60}}
\put(175,-8){\makebox(0,0)[l]{$q^3\mi q^2\pl1$}}
\put(235,-8){\makebox(0,0)[r]{$q^2\mi1$}}
\put(270,0){\oval(70,20)}
\put(270,0){\makebox(0,0){$(q^3\mi q^2\pl1)(q^3\mi1)$}}
\put(270,-15){\makebox(0,0){$q^3\mi q^2$}}
\end{picture}
$$
\vskip 1.0mm}

\Proof
Note that the adjacency relation is symmetric, so that $Z$
is an undirected graph. The computation of the parameters is completely
straightforward. Clearly, $Z$ has $q^6$ vertices.
For $a,b \in W$ the maps $(u,u') \mapsto (u+a,u'+(a \times u)+b)$ are
automorphisms of $Z$, so ${\rm Aut} (Z)$ is vertex-transitive.

The $q^3-1$ neighbours of $(0,0)$ are the vertices $(v,0)$ with $v \ne 0$.
The common neighbours of $(0,0)$ and $(v,0)$ are the vertices $(cv,0)$
for $c \in \Ff_q$, $c \ne 0,1$. Hence $a_1 = q-2$.

The $(q^3-1)(q^2-1)$ vertices at distance 2 from $(0,0)$ are the
vertices $(u,u')$ with $u,u' \ne 0$ and $u' \perp u$.
The common neighbours of $(0,0)$ and $(u,u')$ are the $(v,0)$ with
$v \times u = u'$, and together with $(v,0)$ also $(v+cu,0)$ is a
common neighbour, so $c_2 = q$.
Vertices $(u,u')$ and $(v,v')$, both at distance 2 from $(0,0)$
are adjacent when $0 \ne v \perp u'$ and $v \ne u$ and $v \times u \ne u'$
and $v' = u \times v + u'$, so that $a_2 = q^2-q-2$.

The remaining $(q^3-1)(q^3-q^2+1)$ vertices have distance 3 to $(0,0)$.
They are the $(w,w')$ with $w \not\perp w'$ or $w = 0 \ne w'$.
The neighbours $(u,u')$ of $(w,w')$ that lie at distance 2 to $(0,0)$
satisfy $0 \ne u \perp w'$ and $(0 \ne) u' = w \times u + w'$,
so that $c_3 = q^2 - 1$. This shows that $Z$ is distance-regular
with the claimed parameters. The spectrum follows.

The fact that the extended bipartite double is distance-regular,
and has the stated intersection array, follows from \cite{BCN},
Theorem 1.11.2(vi).

The fact that $Z_3$ is strongly regular follows from
\cite{BCN}, Proposition 4.2.17(ii) (which says that this happens
when $Z$ has eigenvalue $-1$).
\qed

\medskip
For $q=2$, the graphs here are (i) the folded 7-cube,
(ii) the folded 8-cube, (iii) the halved folded 8-cube.
All are distance-transitive.
For $q > 2$ these graphs are not distance-transitive.

When $q$ is a power of two, the graphs $\hat{Z}$ have the same
parameters as certain Kasami graphs, but for $q > 2$ these are nonisomorphic.

\medskip
Next, a more geometric description of this graph.

\medskip
Let $H$ be a vector space of dimension 7 over the field $\Ff_q$,
provided with a nondegenerate quadratic form.
Let $\Gamma$ be the graph of which the vertices are the maximal totally
isotropic subspaces of $H$ (of dimension 3), where two vertices are
adjacent when their intersection has dimension 2.
This graph is known as the dual polar graph of type $B_3(q)$.
It is distance-regular with intersection array
$\{q(q^2+q+1), q^2(q+1), q^3; 1, q+1, q^2+q+1 \}$.
(See \cite{BCN}, \S9.4.)

\begin{Proposition}\label{B3prop}
Let $\Gamma$ be the dual polar graph of type $B_3(q)$.
Fix a vertex $\pi_0$ of $\Gamma$, and let $\Delta$ be the subgraph
of $\Gamma$ induced on the collection of vertices disjoint from $\pi_0$.
Then $\Delta$ is isomorphic to the graph $Z$ of Proposition \ref{B3expl}.
Its extended bipartite double $\hat{\Delta}$ (or $\hat{Z}$) is isomorphic
to the graph of Proposition \ref{D4prop}.
\end{Proposition}

\Proof
Let $V$ be a vector space of dimension 8 over $\Ff_q$
(with basis $\{e_1,\ldots,e_8\}$),
provided with the nondegenerate quadratic form
$Q(x) = x_1x_5+x_2x_6+x_3x_7+x_4x_8$. The point $P = (0,0,0,1,0,0,0,-1)$
is nonisotropic, and $P^\perp$ is the hyperplane $H$ defined by $x_4=x_8$.
Restricted to $H$ the quadratic form becomes
$Q(x) = x_1x_5+x_2x_6+x_3x_7+x_4^2$.

The $D_4$-geometry on $V$ has disjoint maximal totally isotropic
subspaces $E = \langle e_1,e_2,e_3,e_4\rangle$
and $F = \langle e_5,e_6,e_7,e_8\rangle$.
Fix $E$ and consider the collection of all maximal totally isotropic
subspaces disjoint from $E$. This is precisely the collection of
images $F_A$ of $F$ under matrices
$\left(\begin{array}{cc}I&A\\0&I\end{array}\right)$,
where $A$ is alternating with zero diagonal (cf.~\cite{BCN},
Proposition 9.5.1(i)).
Hence, we can label the $q^6$ vertices $F_A \cap H$ of $\Delta$
with the $q^6$ matrices $A$.

Two vertices are adjacent when they have a line in common,
that is, when they are the intersections with $H$ of maximal totally isotropic
subspaces in $V$, disjoint from $E$, that meet in a line contained in $H$.
Let $$A = \left( \begin{array}{cccc}
0 & a & b & c \\
-a & 0 & d & e \\
-b & -d & 0 & f \\
-c & -e & -f & 0
\end{array}\right).$$
Then $\det A = (af-be+cd)^2$, and if $\det A = 0$ but $A \ne 0$, then $\ker A$
has dimension 2, and is spanned by the four vectors
$(0,f,-e,d)^\top$, $(-f,0,c,-b)^\top$,
$(e,-c,0,a)^\top$, $(-d,b,-a,0)^\top$.
Writing the condition that matrices $A$ and $A'$ belong to adjacent vertices
we find the description of Proposition \ref{B3expl} if we take
$u = (c,e,f)$ and $u' = (-d,b,-a)$.
\qed

\section{History}
In 1991 the second author constructed the graphs from Section \ref{bip}
and the first author those from Section \ref{prim}. Both were mentioned
on the web page \cite{ac}, but not published thus far.
These graphs have been called the Pasechnik graphs and the
Brouwer-Pasechnik graphs, respectively, by on-line servers.

\medskip\noindent
Addresses of authors:

\medskip
\begin{minipage}{4.2cm}
Andries E. Brouwer \\
Dept. of Math. \\
Techn. Univ. Eindhoven \\
P. O. Box 513 \\
5600MB Eindhoven \\
Netherlands \\
{\tt aeb@cwi.nl}
\end{minipage}
\begin{minipage}{7.3cm}
Dmitrii V. Pasechnik \\
School of Physical and Mathematical Sciences \\
Nanyang Technological University \\
21 Nanyang Link \\
Singapore 637371 \\
{\tt dima@ntu.edu.sg}
\end{minipage}

\end{document}